\documentclass[10pt]{amsart}

\usepackage[all, cmtip]{xy}
\usepackage{amsmath,amsthm,amssymb,amsfonts,amscd}
\usepackage{tikz-cd}
\usepackage{float} 

\usepackage{graphicx} 
 
\usepackage{amsmath,amscd}
\usepackage[all,cmtip]{xy}

\newtheorem{theorem}{Theorem}[section]
\newtheorem{lemma}[theorem]{Lemma}
\newtheorem{proposition}[theorem]{Proposition}
\newtheorem{corollary}[theorem]{Corollary}

\theoremstyle{definition}
\newtheorem{definition}[theorem]{Definition}
\newtheorem{notation}[theorem]{Notation}

\newtheorem{question}[theorem]{Question}
 \newtheorem{algorithm}[theorem]{Algorithm}

\theoremstyle{remark}
\newtheorem{remark}[theorem]{Remark}

\numberwithin{equation}{section}

\begin{document}

\title[Cascades of toric log del Pezzo surfaces of Picard number one]{Cascades of toric log del Pezzo surfaces \\ of Picard number one}

\author{DongSeon Hwang}
 
\address{Department of Mathematics, Ajou University, Suwon 16499, Republic Of Korea}
 
\email{dshwang@ajou.ac.kr}
 
\thanks{This research was supported by Samsung Science and Technology Foundation under Project SSTF-BA1602-03.}

 
\date{Dec 5, 2020} 


\begin{abstract}
We classify toric log del Pezzo surfaces of Picard number one by introducing the notion,  cascades.     As an application, we show that if such a surface is K\"ahler-Einstein, then it should admit a special cascade, and it satisfies the equality of the orbifold Bogomolov-Miyaoka-Yau inequality, i.e., $K^2 = 3e_{orb}.$  
\end{abstract}

\maketitle

\section{Introduction}
We work over an algebraically closed field of characteristic zero.

A lot of work has been devoted to classify log del Pezzo surfaces. For example, they are classified up to index $3$. See  \cite{HW}, \cite{AN}, \cite{KK}, \cite{Nakayama}, and \cite{FY}. For the Picard number one case, see also \cite{Zhang} and \cite{KM}.  Moreover, toric log del Pezzo surfaces are completely classified up to index $17$ (\cite{KKN}). See   \cite{GRDB} for the list. Recently, all toric log del Pezzo surfaces with $1$ singular point are completely classified in  \cite{Dais}, and those with $2$ are completely classified in \cite{Suyama} which also contains a  partial classification of those with  $3$ singular points.

In this note, we shall classify toric log del Pezzo surfaces of Picard number one by using the notion of a cascade, which was introduced in \cite{H} for a larger class of rational $\mathbb{Q}$-homology projective planes. In fact, even though there exists infinitely many toric log del Pezzo surfaces of Picard number one, one might think that classifying them is not a very difficult task at least  in the sense that it is easy to describe their corresponding Fano triangles, see e.g., \cite[Proof of Proposition 3.10]{HKim}. But it does not give us any geometric intuition and thus sometimes it is not easy to derive geometric consequences. By describing the classification in terms of cascades, which we will soon define, one can understand the underlying geometry more clear. See Section 4 for the applications.  For example, one can easily determine whether a toric log del Pezzo surface of Picard number one with given singularity types exists or not. See Theorem \ref{fiber}, Corollary \ref{classification} and Algorithm \ref{Algorithm}.

\begin{definition}
Let $S$ be a toric log del Pezzo surface of Picard number one. We say that  $S$ \emph{admits a cascade}  if there exists a diagram as follows:

 $$\begin{CD}
S' = S'_t  @>\phi_t >> S'_{t-1}  @>\phi_{t-1} >> \ldots      @>\phi_1>> S'_{0}    \\ 
@V\pi_tVV @V\pi_{t-1}VV @.     @V\pi_0VV \\
S_t   := S @. S_{t-1} @. \ldots     @. S_{0}
\end{CD}$$
 
where for each $k$ 
\begin{enumerate}
\item $\phi_{k}$ is a toric blow-down,  
\item $\pi_k$ is the minimal resolution, 
\item  $S_k$ is a toric log del Pezzo surface of Picard number one,  and 
\item $S_0$ is basic. (See Definition \ref{basic} for the definition.)
\end{enumerate}  
In this case, we also say that $S$ \emph{admits a cascade to $S_0$}, and $S_0$ is the \emph{basic surface} of $S$. 
\end{definition}
 
In the above definition,  when $S_k$ is already smooth, by removing the condition of being Picard number one, we can simply take $S'_k := S_k$ and set  $\pi_k$  to be  the identity morphism. This reminds us the classical construction of smooth del Pezzo surfaces.

The first main result of the present paper is to show the existence of a cascade for every toric log del Pezzo surface of Picard number one.

\begin{theorem}\label{existence}
Every toric log del Pezzo surface of Picard number one admits a cascade. 
\end{theorem}

The proof uses the standard theory of $\mathbb{P}^1$-fibrations. By looking at the dual graph of the torus-invariant divisors,  one can immediately read off the information of $\mathbb{P}^1$-fibration structure on the corresponding smooth toric surface. See Notation \ref{basicdualgraph} for dual graphs.

Conversely, by inverting the cascade process, one can obtain every  toric log del Pezzo surface of Picard number one.

\begin{theorem}\label{inverting}
The minimal resolution of every toric log del Pezzo surface of Picard number one that is not basic is obtained from one of the three basic toric surfaces $S(std^1_n)$, $S(std^2_n)$ and $S(3A_2)$  by a sequence of toric blowups at the intersection point of a  $(-1)$-curve and a torus-invariant  curve with self-intersection number at most $-2$. 
(See Notation \ref{basicsurface}   for the definition of $S(std^1_n)$, $S(std^2_n)$ and $S(3A_2)$.)
\end{theorem}

Since the cascade and its inverse process preserve the number of singular points of $S$, we can describe all   toric log del Pezzo surfaces of Picard number one with respect to the  given number of singular points.

\begin{theorem}
Let $S$ be a toric log del Pezzo surface of Picard number one. If $S$ is not basic,  
it admits a cascade to one of the three basic surfaces: $S(3A_2),$ $S(std^1_n)$, $S(std^2_n)$.  
In particular, we have the following.
\begin{enumerate}
    \item If $|Sing(S)|=0$, then $S \cong \mathbb{P}^2$.
    \item If $|Sing(S)|=1$, then $S \cong \mathbb{P}(1,1,n)$ where $n \geq 2$.
    \item If $|Sing(S)|=2$, then $S \cong \mathbb{P}(1,p,q)$ and it admits a cascade to $S(std^1_n)$.    
    \item If $|Sing(S)|=3$, then $S$ admits a cascade to either $S(std^2_n)$ or $S(3A_2)$.  
\end{enumerate}
\end{theorem}

In particular, this reproves the theorem by \cite{Dais} and \cite{Suyama} for the Picard number one case.

As an application, we consider the orbifold Bogomolov-Miyaoka-Yau inequality. The inequality does not hold in general for Fano manifolds or Fano orbifolds. However, Chan and Leung proposed a Miyaoka-Yau type inequality for K\"ahler-Einstein  toric Fano manifolds.

\begin{theorem}\cite[Theorem 1.2]{CL}
Let $X$ be a K\"ahler-Einstein toric Fano manifold of dimension $n$. Then, for any nef class $H$, we have 
$$c^2_1(X)H^{n-2} \leq 3c_2(X)H^{n-2}$$
if either $n=2,3,4$, or each facet of the corresponding dual polytope of the Fano polytope of $X$ contains a lattice point in its interior.
\end{theorem} 

It is natural to ask whether the above inequality can be generalized in singular setting. 

\begin{question}\cite[Question 1.8]{HKim}\label{KEMY}
Let $S$ be a K\"ahler-Einstein toric log del Pezzo surface. Does the inequality 
$K^2(S) \leq 3e_{orb}(S)$
 holds? 
\end{question}

Unfortunately, the answer  is negative in general as in \cite[Example 1.9]{HKim}. But  it holds  when the Picard number is one. 
 
\begin{theorem}\label{KEtriangleMY} Let $S$ be a K\"ahler-Einstein toric log del Pezzo surface of Picard number one. Then, we have the following properties.
\begin{enumerate}
    \item $K^2_S=3e_{orb}$. 
    \item $S$ is either isomorphic to $\mathbb{P}^2$ or $S$ has exactly $3$ singular points.
    \item If $S$ is not isomorphic to  $\mathbb{P}^2$,  it admits a cascade to $S(3A_2)$, not to $S(std^2_n)$.
\end{enumerate}
\end{theorem}

We emphasize that the condition of being K\"ahler-Einstein forces  a singular toric log del Pezzo surface of Picard number one to admits a cascade to a particular basic surface, i.e., $S(3A_2)$.

As a final application, we give a simple observation that every finite cyclic group is a Brauer group of a toric log del Pezzo surface of Picard number one. See Theorem \ref{Brauerexample}.

\section{Basic toric log del Pezzo surfaces of Picard number one} 
\label{p1fibration}

Throughout this section, we always denote by  $S$  a toric log del Pezzo surface of Picard number one, $f:S' \rightarrow S$ be its minimal resolution. Note that if $S$ is singular, i.e., the Picard number of $S'$ is greater than one,  the torus-invariant divisors form   two sections and two fibers of a suitable $\mathbb{P}^1$-fibration $\Phi: S' \rightarrow \mathbb{P}^1$. For generalities about $\mathbb{P}^1$-fibrations on rational surfaces, see \cite{Miyanishi} or \cite{GMM}.

\begin{notation} 
We denote by $[[s^2_1, F_1, s^2_2, F_2]]$ the smooth toric surface $S'$  admitting a $\mathbb{P}^1$-fibration $\Phi: S' \rightarrow \mathbb{P}^1$ where $s_1$ and $s_2$ are the two torus-invariant sections of $\Phi$ and; $F_1$ and $F_2$  are the two torus-invariant  fibers of $\Phi$.
\end{notation}

\begin{definition}
Let $F$ be a singular fiber  of a $\mathbb{P}^1$-fibration on $S'$. 
    \begin{enumerate}
    \item $F$ is said to be \emph{of type $I_0$} if its dual graph is of the form $\overset{-2}{\circ}-\overset{-1}{\circ}-\overset{-2}{\circ}$.
    \item $F$ is said to be \emph{of type $I$} if it can be contracted to a fiber of type $I_0$.
    \item $F$ is said to be \emph{of type $II_0$} if its dual graph is of the form $\overset{-1}{\circ}-\overset{-2}{\circ}-\overset{-2}{\circ}-\overset{-1}{\circ}$.   
    \item $F$ is said to be \emph{of type $II$} if it can be contracted to a fiber of type $II_0$.
    \end{enumerate}
\end{definition}

\begin{notation} Let $\Phi$ be a $\mathbb{P}^1$-fibration. 
    \begin{enumerate}
    \item A smooth fiber  is denoted by $F_0$.
    \item A singular fiber of type $I_0$ is denoted by $F^0_1$.
    \item A singular fiber of type $I$ is denoted by $F_1$.
    \item A singular fiber of type $II_0$ is denoted by $F^0_2$.
    \item A singular fiber of type $II$ is denoted by $F_2$.     
    \end{enumerate}
\end{notation}

The below lemma immediately follows from the standard theory of smooth projective rational surfaces. 

\begin{lemma}\label{N=12-3n}
Let $S$ be a toric log del Pezzo surface of Picard number one and $S'$ be its minimal resolution. Denote by $n$ the the number of torus-invariant curves on $S'$ and by $N$ the sum of all self-intersection numbers of the torus-invariant curves.  Then, 
we have $N = 12 - 3n.$ 
\end{lemma}

The following notion is essential in the description of the cascades.
\begin{definition} \label{basic}   $S$ is said to be \emph{basic} if $D^2 \geq -2$  for every torus-invariant curve $D$ on $S'$ intersecting $C$ where $C$ is any  $(-1)$-curve on $S'$.  
\end{definition}

For later use,  we introduce the following five surfaces.
 
\begin{notation}\label{basicsurface}
We define the below five surfaces equipped with a $\mathbb{P}^1$-fibration structure.
    \begin{enumerate}
        \item $S(\mathbb{P}^2) = \mathbb{P}^2$
        \item $S(std^0_n) := [[-n,F_0,n,F_0]]$. 
        \item $S(std^1_n) := [[-n, F_0, n-1, F^0_1]]$. 
        \item $S(std^2_n) := [[-n, F^0_1, n-2, F^0_1]]$. 
        \item $S(3A_2) := [[-2, F^0_1, -2, F^0_2]]$. 
    \end{enumerate}   
Figure \ref{basicdualgraph}  describes the \emph{basic dual graphs}, i.e., the dual graphs of the torus-invariant curves on the above five surfaces.



\begin{figure}[H]  \label{basicdualgraph}
\centering
\begin{tikzpicture}[thick,transform canvas={scale=0.5}]
\edef\shiftx{-10.5}; \edef\shifty{-1.8} 

  \begin{scope}[shift={(0+\shiftx,0+\shifty)}]
    \node at (-90: 2) {$G(\mathbb{P}^2)$} ; 
    \tikzstyle{every node} = [draw, shape=circle, fill=white, minimum size=.8cm];
    \node (1) at (30: 1) {$+1$} ;
    \node (2) at (150: 1) {$+1$} ;
    \node (3) at (-90: 1) {$+1$} ; 
    \draw (1) -- (2) -- (3) -- (1) ;
  \end{scope}
%
  \begin{scope}[shift={(5+\shiftx,0+\shifty)}]
    \node at (-90: 2) {$G(std_n^0)$} ; 
    \tikzstyle{every node} = [draw, shape=circle, fill=white, minimum size=.8cm];
    \node (1) at (0: 1) {$0$} ;
    \node (2) at (90: 1) {$-n$} ;
    \node (3) at (180: 1) {$0$} ;
    \node (4) at (-90: 1) {$n$} ;
    \draw (1) -- (2) -- (3) -- (4) -- (1) ;
  \end{scope}
%
  \begin{scope}[shift={(10+\shiftx,0+\shifty)}]
    \node at (-90: {sqrt(2)+1}) {$G(std_n^1)$} ; 
    \tikzstyle{every node} = [draw, shape=circle, fill=white, minimum size=.8cm];
    \node (1) at (0: {sqrt(2)}) {$0$} ;
    \node (2) at (90: {sqrt(2)}) {$-n$} ;
    \node (3) at (135: 2) {$-2$} ;
    \node (4) at (180: {sqrt(2)}) {$-1$} ;
    \node (5) at (-135: 2) {$-2$} ;
    \node (6) at (-90: {sqrt(2)}) {$\scriptstyle n-1$} ;
    \draw (1) -- (2) -- (3) -- (4) -- (5) -- (6) -- (1) ;
  \end{scope}
%
  \begin{scope}[shift={(15+\shiftx,0+\shifty)}]
    \node at (-90: {sqrt(2)+1}) {$G(std_n^2)$} ; 
    \tikzstyle{every node} = [draw, shape=circle, fill=white, minimum size=.8cm];
    \node (1) at (45: 2) {$-2$} ;
    \node (2) at (90: {sqrt(2)}) {$-n$} ;
    \node (3) at (135: 2) {$-2$} ;
    \node (4) at (180: {sqrt(2)}) {$-1$} ;
    \node (5) at (-135: 2) {$-2$} ;
    \node (6) at (-90: {sqrt(2)}) {$\scriptstyle n-2$} ;
    \node (7) at (-45: 2) {$-2$} ;
    \node (8) at (0: {sqrt(2)}) {$-1$} ;
    \draw (1) -- (2) -- (3) -- (4) -- (5) -- (6) -- (7) -- (8) -- (1) ; 
  \end{scope}
%

  \begin{scope}[shift={(20+\shiftx,0+\shifty)}]
    \node at (-90: 3) {$G(3A_2)$} ; 
    \tikzstyle{every node} = [draw, shape=circle, fill=white, minimum size=.8cm];
    \node (1) at (45: {2*sqrt(2)}) {$-1$} ;
    \node (2) at (90: 2) {$-2$} ;
    \node (3) at (135: {2*sqrt(2)}) {$-2$} ;
    \node (4) at (180: 2) {$-1$} ;
    \node (5) at (-135: {2*sqrt(2)}) {$-2$} ;
    \node (6) at (-90: 2) {$-2$} ;
    \node (7) at (-45: {2*sqrt(2)}) {$-1$} ;
    \node (8) at (2,-{2/3}) {$-2$} ;
    \node (9) at (2,{2/3}) {$-2$} ;
    \draw (1) -- (2) -- (3) -- (4) -- (5) -- (6) -- (7) -- (8) -- (9) -- (1) ;
  \end{scope}
 
\end{tikzpicture}
\vspace{2.5cm}
\caption{Basic dual graphs}
\end{figure}
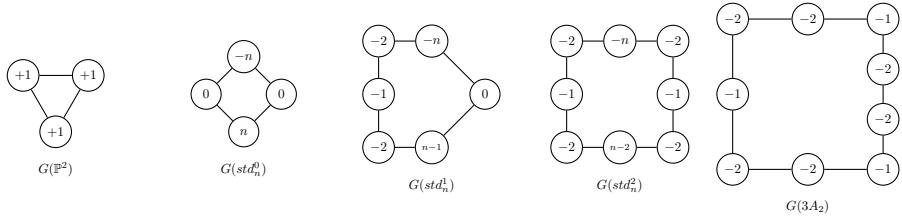
\end{notation}

Now we are ready to determine basic toric log del Pezzo surfaces of Picard number one. 
 
\begin{proposition}\label{typeO}
If $S$ is basic, then its minimal resolution  is  one of the following five surfaces: $S(\mathbb{P}^2)$,  $S(std^0_n)$,      $S(std^1_n)$,    $S(std^2_n) $,    $S(3A_2)$. 
\end{proposition}
 
\begin{proof} 
Assume that $S$ is not isomorphic to  $\mathbb{P}^2$.
Then, $S$ is singular and the Picard number of its minimal resolution $S'$ is greater than one. Thus, $S'$ admits a $\mathbb{P}^1$-fibration $\pi: S' \rightarrow \mathbb{P}^1$ where the cycle of torus-invariant curves forms two singular fibers and two sections of $\pi$.

If $\pi$ is relatively minimal, then $S'$ is isomorphic to the Hirzebruch surface $\mathbb{F}_n = S(std^0_n)$ with $n \neq 1$. In this case, $S$ is isomorphic to $\mathbb{P}(1,1,n)$. 

From now on, we assume that $\pi$ is not relatively minimal. In particular,  there exists a $(-1)$-curve on $S'$. Moreover, since $S$ is singular, there exists a torus-invariant curve with self-intersection number at most $-2$.

Note that  there exists a  $(-1)$-curve $E$ meeting one of the exceptional curves of $f$. Let $D_1, D_2, \ldots, D_k$ be a chain of torus-invariant curves which contracts to one of the singular points of $S$ such that $E$ intersects $D_1$. Let $C$ be the other torus-invariant curve intersecting $E$. Since $S$ is basic, $D^2_1 \geq -2$ and $C^2 \geq -2$. We first consider the case   $D^2_1 = C^2 = -2$. Since $D_1+2E+C$ induces a $\mathbb{P}^1$-fibration structure on $S'$ on which it forms a singular fiber,  there exists another torus-invariant fiber $F$. Since $S$ is basic, it is easy to see that $F$ is one of the following four cases: 
$$\overset{0}{\circ}, \overset{-1}{\circ}-\overset{-1}{\circ}, \overset{-2}{\circ}-\overset{-1}{\circ}-\overset{-2}{\circ}, \overset{-1}{\circ}-\overset{-2}{\circ}-\overset{-2}{\circ}-\ldots-\overset{-2}{\circ}-\overset{-1}{\circ}.$$
In the first and third case, by Lemma \ref{N=12-3n}, we see that the corresponding $\mathbb{P}^1$-fibration structures   are  $S(std^1_n)$ and $S(std^2_n) $, respectively,   where $n \geq 2$. In the second case, one can show that $S$ is of Picard number $2$ or $3$, a contradiction. In the final case, since $S$ is basic and $\rho(S)=1$, the torus-invariant sections have self-intersection number $-2$, so $S$ has only rational double points as singular points. Thus, by Lemma \ref{N=12-3n}, one can see that $F$ should be of type ($II_0$). 
Hence, the corresponding surface is $S(3A_2)$.

Now we consider the case  $C^2 = -1$ and assume that there is no  $(-1)$-curve such that its adjacent torus-invariant curves have self-intersection number at most $-2$.  Then,  since $E+C$ induces a $\mathbb{P}^1$-fibration on which it forms a  complete fiber, there exists another torus-invariant fiber $F$. By assumption, we see that  the fiber $F$ is one of the following:
 $$\overset{0}{\circ}, \overset{-1}{\circ}-\overset{-1}{\circ}, \overset{-1}{\circ}-\overset{-2}{\circ}-\overset{-2}{\circ}-\ldots-\overset{-2}{\circ}-\overset{-1}{\circ}.$$
One can see that $\rho(S) > 1 $ in all of the above cases, which is a contradiction.

Finally, we may assume that, for every  $(-1)$-curve $E$ intersecting an exceptional curve $D_1$ of $f$,  the other torus-invariant curve $C$ intersecting $E$ have $C^2 \geq 0$.  Then, by the similar analysis as above, one can see that $C$ is a section of a $\mathbb{P}^1$-fibration $\Phi$, $D_1$ is part of a fiber of type $II$ and the other fiber is either of type $II$ or of the form  $\overset{0}{\circ}$. In any case, we have   $\rho(S) > 1$,   a contradiction. 
\end{proof}

\begin{remark}
Proposition \ref{typeO} shows that $D^2 \geq -2$ can be replaced by $D^2 = -2$ in Definition \ref{basic}.  
\end{remark}

Every toric log del Pezzo surface of Picard number one corresponds to a Fano triangle. See \cite{KN} for a general introduction to Fano polytopes. For each basic surface $S(X)$  in Notation \ref{basicsurface}, we denote by $P(X)$ the corresponding Fano triangle. See Figure \ref{Polytopecoord} for the    explicit coordinates for the ray generators of $P(X)$ that is basic.

\begin{figure}[H]\label{Polytopecoord}
$
\begin{array}{|r|l|}
\hline
P(\mathbb{P}^2) &  \{(0,1),(1,0),(-1,-1)\}  \\
P(std^0_n), n\geq 2 &  \{ (0,1), (-1,0), (1, -n) \}   \\
P(std^1_n), n\geq 2 & \{ (0,1), (-2,1), (1, -n) \}  \\
P(std^2_n), n\geq2 &  \{ (0,1), (-2, 1), (2, -2n+1) \}  \\
P(3A_2) & \{ (-2, 1), (1, -2), (1,1) \} \\
\hline
\end{array}
$
\caption{Ray generators for $P(X)$ that is basic}
\end{figure}
\begin{figure}[H]\label{RefPolytope}
$\begin{array}{cccc}
P(std^0_2) & P(std^1_2)  &  P(std^2_2) & P(3A_2) \\
\includegraphics[width=2.1cm]{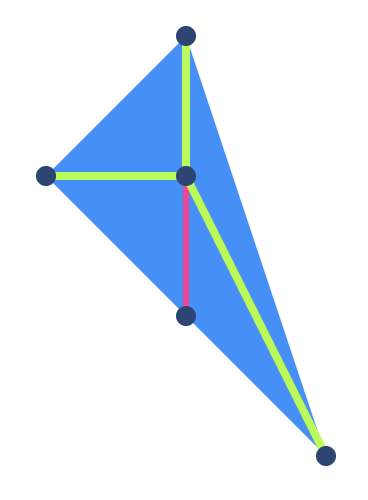} &
\includegraphics[width=3cm]{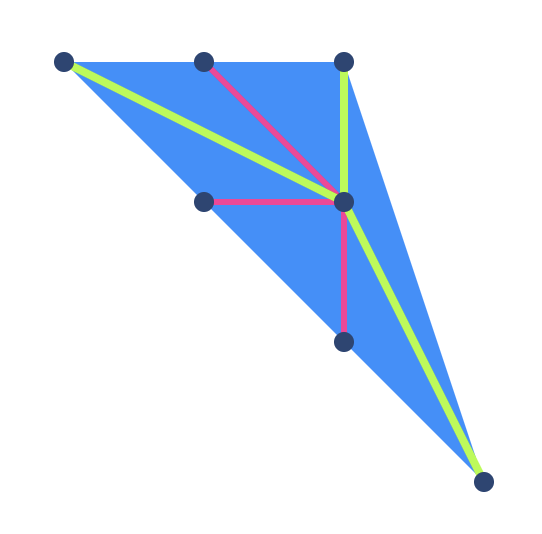} &
\includegraphics[width=3cm]{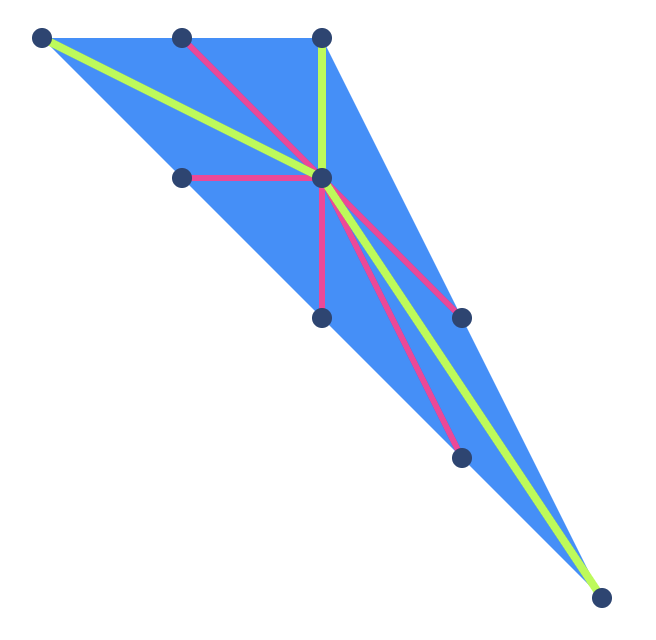} &
\includegraphics[width=3cm]{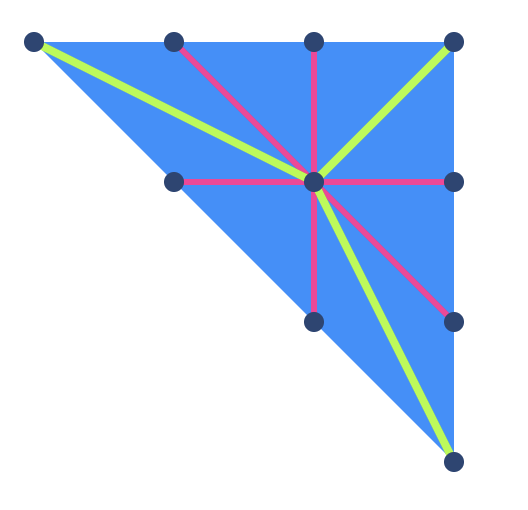}\\
\end{array}$
\caption{Reflexive singular basic Fano polygons}
\end{figure}

See Figure \ref{RefPolytope} for the drawings of  reflexive singular basic Fano triangles.
 
\section{Cascades of toric log del Pezzo surfaces of Picard number one}

\begin{definition}\label{one-step}
Let $S$ be a toric log del Pezzo surface of Picard number one. We say that $S$ \emph{admits a one-step cascade}  if there exists a diagram as follows:

 $$\begin{CD}
S'   @>\phi >> \bar{S'}    \\ 
@V\pi VV @V \bar{\pi} VV      \\
S    @. \bar{S}  
\end{CD}$$
 
where 
\begin{enumerate}
\item $\phi$ is a blow-down of a  $(-1)$-curve, 
\item $\pi$ and $\bar{\pi}$ are minimal resolutions, and 
\item  $\bar{S}$ is a toric log del Pezzo surface of Picard number one.
\end{enumerate}  
\end{definition}

\subsection{Existence of a cascade(=Proof of Theorem \ref{existence})}
If $S$ is basic, we are done. Assume that $S$ is not basic. Then, there exists a  $(-1)$-curve $E$ that intersects a torus-invariant curve $C$ with $C^2 \leq -3$. Let $D$ be the other torus-invariant curve intersecting $E$. We claim that $D^2 = -2$. By \cite[Lemma 1.4]{Zhang},  $D^2 \geq -2$. If $D^2 \geq -1$, then, by contracting $E$ and then contracting all torus-invariant curves with self-intersection number at most $-2$, we get a projective surface of Picard number zero, which is a contradiction.   
Thus, we have $D^2 = -2.$ This can also be derived from \cite[Lemma 4.2]{Zhang}. Now, contracting $E$ induces a one-step cascade. 
  
\subsection{Inverting a cascade(=Proof of Theorem \ref{inverting})}

In the process of each one-step cascade $\phi$, the blowing-up locus of $\phi$, in the notation of Definition \ref{one-step},  is exactly  the intersection point of two  torus-invariant curves, one of them being  a  $(-1)$-curve and the other one has self-intersection number at most $-2$.  
Since there are   exactly three basic surfaces $S(std^1_n)$, $S(std^2_n)$ and $S(3A_2)$ containing a torus-invariant $(-1)$-curve, the result follows.

\subsection{Properties of cascades}   To describe the applications of Theorem \ref{existence} and Theorem \ref{inverting}, we introduce the notion of a trace of $S$.

\begin{definition}
The sum of self-intersection numbers of all irreducible components of exceptional curves of $f$ multiplied by $-1$ is called the \emph{trace} $tr(S)$ of $S$. In other words,
$$ tr(s) = -\sum D^2_i$$
where the sum runs over all exceptional curves $D_i$ over $f$.
\end{definition}

Now, the proof of Theorem \ref{inverting} immediately yields the following. 

\begin{corollary}\label{one-stepprop} The the number of singular points of $S$, and the number $tr(S)-3L$  are invariant under a  cascade where $L$ denotes the number of exceptional curves of the minimal resolution.  
\end{corollary}

\begin{proof}
It is enough to consider only a one-step cascade. 
It is clear that $tr(S)-3L$ is invariant under a one-step cascade. This also follows from   Lemma \ref{N=12-3n}.  
Note that a one-step cascade does not increase the  number of singular points. Assume that a one-step cascade decrease the number of singular points of $S$. Then, there is a chain of torus-invariant curves whose dual graph is of the form
$$\overset{-n}{\circ}-\overset{-1}{\circ}-\overset{-2}{\circ}-\overset{-m}{\circ} $$
where $m \geq -1$ and $n \leq -3$   since $S$ is not basic. Let $E$ be the $(-1)$-curve in the dual graph intersecting the $(-n)$-curve. By blowing down $E$, we can see that $n = 3$ by \cite[Lemma 4.2]{Zhang}, hence we get the following dual graph.
$$\overset{-2}{\circ}-\overset{-1}{\circ}-\overset{-m}{\circ} $$
for some $m \geq -1$ (cf. \cite[Lemma 1.4]{Zhang}). This cannot be possible since the Picard number is one.
\end{proof}

\begin{remark}\label{trace}
By Corollary \ref{one-stepprop}, we can easily compute the trace of toric log del Pezzo surface of Picard number one once we know its basic surface.

     $$\begin{array}{|c||c|c|c|c|c|}
    \hline
    \text{S} & S(\mathbb{P}^2) & S(std^0_n) & S(std^1_n) & S(std^2_n) & S(3A_2)\\
    \hline
    \text{tr(S)} &  -3 & n & 3L-5+n (\geq 3L-3)  & 3L-7+n (\geq 3L-5) & 3L-6 \\
    \hline
    \end{array}$$

The above table shows that the number of singular points and the trace of $S$ determines uniquely the original surface $S$ and its basic surface and vice versa. See Algorithm \ref{Algorithm}.
\end{remark}

\section{Applications}
We completely classify toric log del Pezzo surfaces of Picard number one and their dual graphs. 

\subsection{Classification}  
 
\begin{theorem} \label{fiber}
Let $S$ be a toric log del Pezzo surface of Picard number one. Then,  
    \begin{enumerate} 
    \item Either $S \cong \mathbb{P}(1,1,n)$ with $n \geq 1$, or $S$ admits a cascade to one of the following:  $S(std^1_n)$, $S(std^2_n)$, and $S(3A_2)$ where $n \geq 2$. 
    \item Let $T$ be the basic surface of $S$. Then, we have the following.
        \begin{enumerate}
            \item If $T=S(std^1_n)$, then $S'=[[-n, F_1, n-1, F_0]]$.
            \item If $T=S(std^2_n)$, then $S'=[[-n,F_1,n-2,F'_1]]$. 
            \item If $T=S(3A_2)$, then $S'=[[-n,F_1, -m, F_2]]$.
        \end{enumerate}
    where   $F'_1$ is a fiber of type $I$ and $F_i$ is a fiber of type $I$ for each $i$.
    \end{enumerate}
\end{theorem}

\begin{proof} 
Since $S(std^0_n) \cong \mathbb{P}(1,1,n),$ (1) immediately follows from Theorem \ref{existence} and Theorem \ref{inverting}. 

We may assume that $S$ is not basic. Then, by taking a finite number of one-step cascade,  we can always find three torus-invariant curves whose dual graph is of the form $\overset{-2}{\circ}-\overset{-1}{\circ}-\overset{-2}{\circ}$. Note that they induce a $\mathbb{P}^1$-fibration $\Phi$ on the minimal resolution $S'$ of $S$, on which they form a singular fiber $F$ of type ($I_0$).

Consider the case  $T = S(std^1_n)$.  Since the inverting process only changes the singular fiber $F$ of the $\mathbb{P}^1$-fibration, we see that  $G(S) = [[-n,F_1,n-1, 0]]$ for some integer $n \geq 2$ with the unique singular fiber   $F_1$ of type $I$.

Consider the case  $T = S(std^2_n)$. Then,  only the two torus-invariant sections of $\Phi$ are invariant under the inverse process among all torus-invariant curves.  Thus, we have $G(S)=[[-n,F_1,n-2,F'_1]]$ for some integer $n \geq 2$ where both $F_1$ and $F'_1$ are  of type $I$.

Consider the case  $T = S(3A_2)$. Since no torus-invariant curve is invariant under the process of cascades in general, the result follows.
  
\end{proof}

Now we classify toric log del Pezzo surfaces   of Picard number one with given number of singular points.

\begin{corollary}\label{classification}
Let $S$ be a toric log del Pezzo surface of Picard number one. Then, we have the following.
\begin{enumerate} 
\item  If $|Sing(S)| \leq 1$, then $S \cong \mathbb{P}(1,1,n)$ for some $n \geq 1$. 
\item If $|Sing(S)| = 2$, then $S \cong \mathbb{P}(1, q, (n-1)q + q_1)$ where $gcd(q,q_1) = 1$.
\item If $|Sing(S)| = 3$, then $S$ is obtained by inverting a cascade from $S(3A_2)$ or $S(std^2_n)$.
\end{enumerate}
In particular, if $|Sing(S)| \leq 2$, then $S$ is a weighted projective plane. 
\end{corollary}

To prove Corollary \ref{classification}, we recall the Hirzebruch-Jung continued fraction.

\begin{definition} For integers $n_1$, $n_2$, \ldots, $n_l$, we set the following notation, 
$$[n_1, n_2, \ldots, n_l] := n_1 - \dfrac{1}{n_2 - \dfrac{1}{\ddots - \dfrac{1}{n_l}}}.$$
If $n_i \geq 2$ for each $i$, then it is called a \emph{Hirzebruch-Jung continued fraction}.
\end{definition}

\begin{proof}[Proof of Corollary \ref{classification}]
For (1) and (3),  the result follows from Theorem \ref{fiber} and  Corollary \ref{one-stepprop}. Assume that $|Sing(S)|$ $= 2$. By Theorem  \ref{fiber} and  Corollary \ref{one-stepprop}, $S$ is obtained by inverting the cascade  from $S(std^1_n)$. 

Let $F$ be a singular fiber of a $\mathbb{P}^1$-fibration of the form 
$$\underset{F_1}{\overset{}{\Box}}-\underset{-1}{\circ}-\underset{F_2}{\overset{}{\Box}}$$
where $F_1$ is the dual graph corresponding to the Hirzebruch-Jung continuned fraction $[n_1, \ldots, n_l]$ and    $F_2$ corresponds to  $[m_1, \ldots, m_t]$.
By Lemma \ref{cf1} below, we let $\frac{Q}{q} = [n, n_1, \ldots, n_l]$ and $\frac{q}{q_1} =  [m_1, \ldots, m_t]$  such  that  $[n_1, \ldots, n_l, 1, m_1, \ldots, m_t] = 0$.

We want to show that $S \cong \mathbb{P}(1, q, (n-1)q+q_1)$. Again, by Lemma \ref{cf1}, we see that 
$$Q = [n,n_1, \ldots, n_l] = n - \frac{q-q_1}{q} = \frac{(n-1)q+q_1}{q}.$$ Thus, $S$ and $\mathbb{P}(1,q,Q)$ have the same singularity types. This completes the proof since the singularity type uniquely determines the surface when $|Sing(S)|=2$.  
\end{proof}

\begin{remark}
It is well known that a weighted projective plane is a toric log del Pezzo surface of Picard number one. One can easily construct infinitely many toric surfaces of Picard number one which is not a weighted projective plane by inverting the cascade  from $S(3A_2)$. See the construction in the proof of Theorem \ref{Brauerexample}.
\end{remark}

\begin{remark}
Corollary \ref{classification} reproves the results in \cite{Dais} and \cite{Suyama} for the Picard number one case. 
\end{remark}

\begin{lemma}\label{cf1}
 Let $[n_1, \ldots, n_l]$ and $[m_1, \ldots, m_t]$ be Hirzebrugh-Jung continued fractions such that $[n_1, \ldots, n_l, 1, m_1, \ldots, m_t]=0$. If 
 $[m_1, \ldots, m_t] = \frac{q}{q_1}$, then $[n_l, \ldots, n_1]$ $= \frac{q}{q-q_1}.$
\end{lemma}

\begin{proof}
This lemma is   well-known and easy to prove. See \cite[Example 1]{R} for the algorithm to compute $[n_1, \ldots, n_l]$  for a given  $[m_1, \ldots, m_t]$.  
\end{proof}

\begin{algorithm}\label{Algorithm}
By Corollary \ref{classification}, we can determine whether there exists a toric log del Pezzo surface $S$ of Picard number one with given singularity types.\\

\underline{INPUT}: an $k$-tuple of rational numbers $(\frac{n_1}{m_1}, \ldots, \frac{n_k}{m_k}) $ where $k$ denotes the number of singular points of $S$ and each rational number describes the singularity type. \\

\underline{OUTPUT}: False if there exists no toric log del Pezzo surface of Picard number one having the given singularity type in INPUT. If it exists, we return $S$ if $S$ is basic, or $S$ and its basic surface  if otherwise.\\

\underline{PROCEDURE}: (using notation in Remark \ref{trace})
\begin{enumerate} 
    \item  If the input is empty, i.e., $k=0$, then $S = \mathbb{P}^2$.
    \item If $k=1$ and $m_1 = 1$, then $S = \mathbb{P}(1,1,n_1)$.
    \item If $k \geq 2$, then reorder the $k$-tuple so that $i \geq j$ if and only if either $n_i > n_j$, or $n_i=n_j$ and $m_i \geq m_j$.
    \item If $k=2$, $m_1 = n_2$ and $\frac{n_1-m_1}{m_1}$ is a positive integer, then $S  \cong \mathbb{P}(1, m_1, n_1)$.
    \item If $k=3$ and $tr=3L-6$, then consider the three dual graphs  of the singularities  corresponding to the triple in INPUT.  Form a cycle $G$ by adding one vertex of weight $-1$ between any two of the three dual graphs. Note that there are four possible ways for forming the cycle.  If the graph is $G(3A_2)$ after a finite number of "blowing-down" of the graph, then $S$ is the toric log del Pezzo surface of Picard number one whose dual graph of the torus-invariant divisors is  $G$.  
    \item If $k=3$ and $tr \geq 3L-5$, then consider the three dual graphs  $G_1$, $G_2$, $G_3$ of the singularities  corresponding to $\frac{n_1}{m_1}, \frac{n_2}{m_2}, \frac{n_3}{m_3}$.  Form a tree $G$ by adding one vertex of weight $-1$ between $G_1$ and $G_2$; and between $G_1$ and $G_3$. If the graph is $G(std^2_n)$ after a finite number of "blowing-down" of the graph, then $S$ is the toric log del Pezzo surface of Picard number one whose dual graph of the torus-invariant divisors is  $G$.  
    \item Return False.
\end{enumerate}
\end{algorithm}

\subsection{K\"ahler-Einstein toric log del Pezzo surfaces of Picard number one}

\begin{proof}[Proof of  Theorem \ref{KEtriangleMY}] Since $\mathbb{P}^2$ is K\"ahler-Einstein, (2) follows from \cite[Corollary 3.11]{HKim}.

Let $S$ be a K\"ahler-Einstein log del Pezzo surface of Picard number one. It is enough to assume that $S$ is singular. Consider the minimal resolution $f:S'\rightarrow S$ of $S$. 
Let  $D_1,  \ldots,  D_L$ be the all  irreducible components of the reduced part $\mathcal{D}$ of the $f$-exceptional divisor. Then, by \cite[Remark 3.12]{HKim}, $S$ has $3$ singular points, each of which has local fundamental group of order $a$. 
Then, by \cite[Section 3 and Lemma 3.6]{HK2011a}, 
$$ K^2_S = tr -3L + 6 + \frac{9}{a}$$ 
where $ tr = -\overset{L}{\underset{k=1}{\sum}} D^2_k$.  
Since $3e_{orb}=\frac{9}{a}$, we see that $K^2_S = 3e_{orb}$ if and only if $tr = 3L-6$ if and only if  $S$ admits a cascade to $S(3A_2)$. The last equivalence follows from Corollary  \ref{classification} and Lemma \ref{trace}. Thus, it remains to show that $S$ admits a cascade to $S(3A_2)$. Now the below lemma completes the proof by Corollary \ref{classification}. 
\end{proof}

\begin{lemma}
Let $S$ be a log del Pezzo surface of Picard number one. If $S$ admits a cascade to $S(std^2_n)$, then $S$ is not K\"ahler-Einstein.
\end{lemma}

\begin{proof}
Let $P$ be the Fano polygon corresponding to $S$. It is enough to show that the barycenter of $P$ is not the origin by \cite[Theorem 1.2]{BB}. Since $S$ admits a cascade to $S(std^2_n)$,  $P$ admits a cascade to $P(std^2_n)$. 
Note that the barycenter of $P(std^2_n)=conv\{(1,-1),(1,1),(-3,-1)\}$ is $(-\frac{1}{3},-\frac{1}{3})$. Since the $y$-coordinate of the barycenter is not increasing during the inverting process of the cascade, the barycenter of $P$ cannot be the origin.
\end{proof}

\subsection{Brauer groups} The Brauer group of a toric surface can easily be computed  by the following theorem. 

\begin{theorem}\cite[Corollary 2.9]{DF}\label{toricbrauer}
Let $X$ be a toric surface,  $\Delta$ be the corresponding  complete fan on $\mathbb{R}^2$  and $\Delta(1) = \{ \rho_1, \ldots, \rho_n \}$.
If $N' = \langle \rho_1 \cap N, \ldots, \rho_n \cap N \rangle$, then $B(X) \cong N/N'$.
\end{theorem}

Now, as an application of the cascade structure, we show that every finite cyclic group is a Brauer group of a toric log del Pezzo surface of Picard number one.

\begin{theorem}\label{Brauerexample}
For each positive integer $n$, there exists a toric log del Pezzo surface $S$ of Picard number one with $Br(S) \cong \mathbb{Z}/n\mathbb{Z}$.
\end{theorem}

\begin{proof} 
First, we observe that $Br(\mathbb{P}^2)$ is a trivial group and $Br(std^2_2) \cong \mathbb{Z}/2$. For each integer $n \geq 3$,  
we shall explicitly construct a toric log del Pezzo surface $S$ of Picard number one with $Br(S) \cong \mathbb{Z}/n\mathbb{Z}$ by inverting the cascade from $S(3A_2)$. Let $S_0 = S(3A_2)$ and $f:S'_0 \rightarrow S_0$ be its minimal resolution. Choose a chain of two  $(-2)$-curves $C_1$ and $C_2$. Let $E_i$ be a $(-1)$-curve intersecting $C_i$ for each $i=1,2$. Blow up the intersection point of $C_1$ and $E_1$, and then blow up the intersection point of $C_2$ and $E_2$. Let $S'_1$ be resulting surface and $S_1$ be its anticanonical model. Note that there exists a $(-1)$-curve $E'_i$ intersecting the proper transform $C'_i$ of $C_i$ for $i=1,2$. Blow up the intersection point of $C'_1$ and $E'_1$, and then blow up the intersection point of $C'_2$ and $E'_2$. Let $S'_2$ be resulting surface and $S_2$ be its anticanonical model. One can continue this process. Note that $S_n$ is a toric log del Pezzo surface of Picard number one with 3 singular points of type $2A_{n+2} + [n+2, n+2]$.  Now it is easy to see that $Br(S_n) \cong \mathbb{Z}/(n+3)\mathbb{Z}$ by Theorem \ref{toricbrauer}.
\end{proof}

\bigskip {\bf Acknowledgements.}   
The author was supported by Samsung Science and Technology Foundation under Project SSTF-BA1602-03.	 
\bigskip



\begin{thebibliography}{99}


\bibitem[AN]{AN} V. Alexeev and V. V. Nikulin, {\it Del Pezzo and K3 surfaces}, MSJ Memoirs, 15. Mathematical Society of Japan, Tokyo, 2006.
 
\bibitem[BB]{BB} R. J. Berman and B. Berndtsson, \textit{Real Monge-Amp??re equations and K??hler-Ricci solitons on toric log Fano varieties}, Ann. Fac. Sci. Toulouse Math.,  (6) \textbf{22} (2013), no. 4, 649???711. 
  


\bibitem[CL]{CL} K. Chan, N. C. Leung, \textit{Miyaoka-Yau-type inequalities for K\"ahler-Einstein manifolds},
Commun. Anal. Geom. \textbf{15} (2007), 359-379.
 
\bibitem[D]{Dais} D. I. Dais, \textit{Toric log del Pezzo surfaces with one singularity},  Adv. Geom. \textbf{20} (2020), no. 1, 121-138. 
 
\bibitem[DF]{DF} F. R. Demeyer and T. J. Ford,  \textit{On the Brauer group of toric varieties}, Trans. Amer. Math. Soc. \textbf{335} (1993), no. 2, 559-577. 
 
\bibitem[FY]{FY} K. Fujita, K. Yasutake,  \textit{Classification of log del Pezzo surfaces of index three}, J. Math. Soc. Japan 69 (2017), no. 1, 163???225.


\bibitem[GMM]{GMM} R.V. Gurjar,  K. Masuda,  M. Miyanishi, \textit{Affine space fibrations},   Polynomial Rings and Affine Algebraic Geometry, 151-193,  Springer Proceedings in Mathematics \& Statistics, vol 319. Springer, 2020.
 
\bibitem[GRDB]{GRDB} Graded Ring Database, http://www.grdb.co.uk/
 

 
\bibitem[H]{H} D. Hwang, \textit{Algebraic Montgomery-Yang problem and cascade structure}, preprint available upon request.
 
\bibitem[HK]{HK2011a} D. Hwang and J. Keum, \textit{The maximum number of singular points on rational homology projective planes}, J. Algebraic Geom., \textbf{20} (2011), 495???523. 
 

\bibitem[HKi]{HKim} D. Hwang and Y. Kim, \textit{Symmetric and K\"ahler-Einstein  toric log del Pezzo surfaces}, preprint available upon request.
 
\bibitem[HW]{HW} F. Hidaka, K. Watanabe,  \textit{Normal Gorenstein surfaces with ample anti-canonical divisor}, Tokyo J. Math. 4 (1981), no. 2, 319???330.

 
\bibitem[KK]{KK} G. Kapustka and M. Kapustka,   \textit{Equations of log del Pezzo surfaces of index $\leq 2$}, Math. Z. \textbf{261} (2009), no. 1, 169???188.
 

\bibitem[KKN]{KKN}A. M. Kasprzyk,   M. Kreuzer and B. Nill B,  \textit{On the combinatorial classification of toric log del Pezzo surfaces}, LMS Journal of Computation and Mathematics 13 (2010), 33-46.

\bibitem[KM]{KM} S. Keel and J. McKernan, {\it Rational curves on quasi-projective surfaces}, Mem. Amer. Math. Soc. {\bf 140} (1999), no. 669.

\bibitem[KN]{KN} A. M. Kasprzyk,  B. Nill,   \textit{Fano polytopes}, Strings, gauge fields, and the geometry behind, 349-364, World Sci. Publ., Hackensack, NJ, 2013.


\bibitem [M]{Miyanishi} M. Miyanishi, \textit{Open Algebraic Surfaces}, CRM Monograph Series, \textbf{12} American Math. Soc. 2001.

\bibitem [N]{Nakayama} N. Nakayama, \textit{Classification of log del Pezzo surfaces of index two}, J. Math. Sci. Univ. Tokyo 14 (2007), no. 3, 293???498.


\bibitem[R]{R} O. Riemenschneider, \textit{Deformationen von Quotientensingularit\"aten (nach zyklischen Gruppen)}, Math. Ann. \textbf{209} (1974), 211-248. 
 

\bibitem[S]{Suyama} Y. Suyama,  \textit{Classification of toric log del Pezzo surfaces with few singular points}, available at 	arXiv:1910.00206.
 

\bibitem[Z]{Zhang} D. Zhang, \textit{Logarithmic del Pezzo surfaces of rank one with contractible boundaries},  Osaka J. Math. \textbf{25} (1988), no. 2, 461-497.

\end{thebibliography}
\end{document}